\documentclass[english]{amsart}
\usepackage[T1]{fontenc}
\usepackage[latin1]{inputenc}
\usepackage{babel}

\makeatletter

\newcommand{\LyX}{L\kern-.1667em\lower.25em\hbox{Y}\kern-.125emX\spacefactor1000}

\theoremstyle{plain}
\newtheorem{thm}{Theorem}[section]
\numberwithin{equation}{section} 
\numberwithin{figure}{section} 
\theoremstyle{plain}
\newtheorem{lem}[thm]{Lemma} 
\theoremstyle{remark}
\newtheorem{rem}[thm]{Remark}

\makeatother

\begin{document}

\title{Semi-local invariants for non resonant Poisson structures on \protect\( S^{1}\! \! \times \! \bf R^{n}\protect \). }

\author{Olivier BRAHIC\protect\( \qquad \protect \)Jean-Paul DUFOUR.}

\begin{abstract}
The local study of a Poisson structure near a ``generic'' isolated zero has been the
subject of many works. The purpose of this paper is to give a ``semi-local'' study
along ``generic'' closed curves of such zeros: we formally classify Poisson
structures defined in a neighborhood of \( \Gamma \! =S^{1}\! \times \{0\} \) in \(
S^{1}\! \times \mathbb{R}^{n} \), that vanish on \( \Gamma \), and whose linear
approximation at one point \( m \) of \( \Gamma  \) is isomorphic to the dual of a
non-resonant Lie algebra.
\end{abstract}
\maketitle

\section{Introduction.}

Let $\Pi$ be a Poisson structure on some manifold $M$, vanishing at a point $m\in M$. Denote $
{\mathcal D}_{\omega }\Pi$ the modular vector field ({[}W2{]}) with respect to some density
\( \omega  \), and suppose that $
{\mathcal D}_{\omega }\Pi$ doesn't vanish at $m$; this requirement is independent of the
chosen density, and forces the Poisson structure to vanish all along the orbit $\Gamma$  of ${\mathcal D}_{\omega }\Pi$ through \( m \), thus $\Gamma$ is a whole curve of singularities for $\Pi$. 
In this paper, we are interested in describing generic Poisson structure in the neighborhood of such $\Gamma$, when it is supposed to be closed, that is $\Gamma=S^1$. Also we will consider
that \( \Pi  \) is defined on a neighborhood of \( \Gamma =S^{1}\! \! \times \!
\{0\} \) in \( S^{1}\! \! \times \! \mathbb{R}^{n} \) parametrized by \( (\theta
,x_{1},\dots ,x_{n}) \), and that \( \Pi  \) is zero on \( \Gamma  \).

The first invariant attached to this situation is the period \( c \) of this~ orbit
\( \Gamma  \): it doesn't depend on the chosen density because \( \Pi  \) vanishes
on it. Recall that this is one of the invariants used by O. Radko ({[}R{]}) to
classify Poisson structures on surfaces.

The second~ thing we have to do is to take care of the linear part of \( \Pi  \) at
the points of \( \Gamma  \). If we choose local coordinates \(
(\widetilde{x}_{0},\widetilde{x}_{1},\dots ,\widetilde{x}_{n}) \) vanishing at \( m
\), and such that \( {\mathcal D}_{\omega }\Pi  \) is \( \frac{\partial }{\partial
\widetilde{x}_{0}} \), this linearized Poisson structure satisfies relations:
\[
\begin{array}{c}
\{\widetilde{x}_{0},\widetilde{x}_{i}\}=\sum ^{n}_{j=1}a_{i,j}\widetilde{x}_{j}\\
\{\widetilde{x}_{r},\widetilde{x}_{i}\}=\sum ^{n}_{j=1}c_{r,i}^{j}\widetilde{x}_{j},
\end{array}\]
 for \( r \) and \( i \) varying from 1 to \( n \). This means that it corresponds
to a Lie algebra (see {[}W1{]}) which is a semi-direct product of \( \mathbb{R} \)
with a \( n \)-dimensional Lie algebra \( A \). In the sequel we will denote by \( D
\) the derivation of \( A \) under which \( \mathbb{R} \) acts on \( A \): in the
above coordinates, it has the matrix with coefficients \( a_{i,j} \). We see also
that, up to a linear isomorphism, this Lie algebra doesn't depend on the point \( m
\) chosen on \( \Gamma  \) as \( {\mathcal D}_{\omega }\Pi  \) is an infinitesimal isomorphism.
So this isomorphism class of Lie algebra is the second invariant we can attach to
our situation.

The generic condition we will impose on $\Pi$ in the sequel, is that there are no resonance relation between the eigenvalues \(
\lambda _{1},\dots ,\lambda _{n} \) of \( D \); by this we mean precisely that there
are no relations
\[ \lambda_i=
p_{1}\lambda _{1}+\cdots +p_{n}\lambda _{n},\]
 \[ \lambda_i+\lambda_j=
p_{1}\lambda _{1}+\cdots +p_{n}\lambda _{n}\] for every $i$ and $j,$ ($i\neq j),$
where \( (p_{1},\dots ,p_{n}) \) is a multi-index with non-negatives $p_k$ and
 \( \sum_{k=1}^{n} p_{k}\geq 2, \) except trivial  relations
 $\lambda_i+\lambda_j=\lambda_i+\lambda_j.$ In particular, these eigenvalues are {\it distinct}, and $D$ is diagonalizable. It is then easy to verify that, in the complex case, the brackets of such linear Poisson structures, so called \bf non resonant\rm, and denoted \(
(\mathbb{C}^{n+1})_{(\lambda _{1}...\lambda _{n})}^{*} \), can always be written as:
\begin{eqnarray*}
\{e_{i},e_{j}\}&=&0,\: \: \: \: \; \; \ i,j=1...n \\
\{e_{0},e_{i}\}&=&\lambda _{i}e_{i},\: \: \; \,  i=1...n, \\
 &
\end{eqnarray*}
(the derivation \( D \) is just \( \{e_{0},-\} \)) so that \( A \) has to be
commutative.

Note that two such Lie algebras \( (\mathbb{C}^{n+1})_{(\lambda _{1}...\lambda
_{n})}^{*} \) and \( (\mathbb{C}^{n+1})_{(\lambda '_{1}...\lambda '_{n})}^{*} \) are
isomorphic if and only if \( \lambda _{i}=d\lambda '_{i} \) , up to a re-indexation
of the \( \lambda _{i} \) 's, and for \( d\in \mathbb{C} \) (see the proof of lemma
2.1). Also, in terms of classification of Lie algebras, the \( \lambda _{i} \) 's
should be seen as a set defined up to multiplication by a scalar of all its
elements.

In {[}D-Z{]} there is a local study of this situation (in the real and complex
cases): it appears there that the Poisson structure is, at least formally,
``quadratizable'', i.e. we can find coordinates such that our Poisson structure has
only linear and quadratic terms. In this paper we will show that this is also true
in a neighborhood of the whole curve \( \Gamma  \). Moreover we will describe all
the invariants attached to the isomorphism class of the germ of \( \Pi  \) along \(
\Gamma  \). However, in order to simplify the study, we will restrict ourselves to
the case where {\bf the eigenvalues $\lambda_i$ are all real}.  In that case the
Lie algebra has the above form but in real coordinates; it will be then denoted by
\( (\mathbb{R}^{n+1})_{(\lambda '_{1}...\lambda '_{n})}^{*} \).

In the analytic context, we will need an extra assumption on these eigenvalues : we will suppose that they verify Bruno's condition, that is : if $I$ denote the set of all mutli-indexes $(c_i)_{i=1...n}\in {\mathbb Z}^n$ such that $c_i\leq-1, \forall i=1\dots n$ and $\sum_{i=1\dots n} c_i \lambda_i\neq 0$, one puts :
\begin{equation}\label{bruno}
\omega_k=min\{|\sum_{i=1}^n c_i \lambda_i|;(c_i)\in I, \sum_{i=1}^n c_i<2^k\},
\end{equation}
Imposong Bruno's condition is requiring that $\sum_{k=1}^\infty \frac{1}{2}log\frac{1}{\omega_k}<\infty$, it is a standard condition for assuring analytic convergence of formal series. 

 It is willingly that we didn't precise if $\Pi$ was supposed to be smooth or analytic, next section is valid for the smooth case, and also, under Bruno's condition, in the analytic case (of course, flat terms $o_{\infty }(x)$ will have to be omitted in the analytic case).

\section{The normal form.}

In this section, we give several lemmas, that will, step by step,  lead to a normal
form for \( \Pi  \) in a neighborhood of \( \Gamma \), with underlying preoccupation
to obtain coordinates that are global with respect to \( \theta \in S^{1} \).

\begin{lem}
The brackets induced by \( \Pi  \), can, in a neighborhood of \( \Gamma  \), be
written as:
\begin{eqnarray*}
\{x_{i},x_{j}\}&=&o_{2}(x),\qquad \qquad \qquad \, \; \; \; \; \; \, \: \: \: \; \;  i,j=1\dots n \\
\{\theta ,x_{i}\}&=&\sum _{j=1}^{n}h_{i,j}x_{j}+o_{2}(x)\qquad \; \quad \quad  i=1\dots n,  
\end{eqnarray*}
where \( h_{i,j} \) are smooth functions on \( S^{1} \), \( o_{2}(x) \) denote
smooth functions on \( S^{1}\times \mathbb{R}^{n} \) and of order two in the
variables \( x_{1},\dots \, ,x_{n} \) , and the matrix \( H_{\theta }=(h_{i,j}(\theta
))_{i,j=1\dots n} \) has eigenvalues \( \{k(\theta )\lambda _{i}/i=1\dots n\} \) for some
\( k:S^{1}\longmapsto \mathbb{R} \).
\end{lem}
\begin{proof}
\( \Pi  \) vanishes on \( \Gamma =\{x_{i}=0\ ;\ i=1\dots n\} \), so that we can write
the brackets as:
\begin{eqnarray*}
\{x_{i},x_{j}\}&=&\sum _{k=1\dots n}u_{k}^{i,j}(\theta )x_{k}+o_{2}(x),\, \qquad  i,j=1\dots n \\
\{\theta ,x_{i}\}&=&\sum _{j=1\dots n}h_{i,j}(\theta )x_{j}+o_{2}(x),\quad \qquad  i=1\dots n.  \\
 &
\end{eqnarray*}
 In fact the lemma is just a consequence of the fact that two semi-direct products
of Lie algebras \( \mathbb{R}\: _{L}\triangleright \! \! \! <\mathbb{R}^{n} \)and \(
\mathbb{R}\: _{L'}\triangleright \! \! \! <\mathbb{R}^{n} \)are isomorphic if and
only if \( L \) and \( L' \) are conjugated up to a scalar, but let 's see this in
detail.

Let \( \theta _{0}\in S^{1} \) fixed, denote \( \Pi ^{(1)}_{\theta _{0}} \) the
linear approximation of \( \Pi  \) at \( (\theta _{0},0,\dots \, ,0) \), according to
the introduction, it is isomorphic to the dual of a Lie algebra with commutative
derived ideal , so all \( u^{i,j}_{k}(\theta ) \) must be zero.

Now, let \( \Psi _{\theta _{0}}\! \!:\! (\mathbb{R}^{3})^{*}_{(\lambda
_{1}\dots \lambda _{n})}\! \longmapsto \Pi ^{(1)}_{\theta _{0}} \) a Poisson
isomorphism, and $$ \widetilde{\Psi }_{\theta _{0}}=\left( \begin{array}{cccc}
k & v_{1} & \dots  & v_{n}\\
w_{1} &  &  & \\
\vdots &  & \Psi _{i,j} & \\
w_{n} &  &  &
\end{array}\right)  $$ the associated matrix; \( (\Psi _{\theta _{0}})^{*} \) being a Lie algebra
isomorphism, it respects derived ideals, and \( (w_{1},\dots \, ,w_{n}) \) must
be zero.

Then we see that \( \det(\Psi _{\theta _{0}})\! =\! k.\det\left( \psi _{i,j}\right)
_{i,j} \), in particular \( k\neq 0 \) and \( \left( \Psi _{i,j}\right) _{i,j} \) is
invertible.

At that point, the conditions for \( \Psi _{\theta _{0}} \) to be Poisson are:
\[
\sum _{k=1,\dots ,n}h_{i,k}\psi _{k,j}=k.\psi _{i,j}\lambda _{i},\qquad  i,j=1\dots n,\]
 and we see that these relations can be written as:
\[
H_{\theta _{0}}.\left( \psi _{i,j}\right) =\left( \psi _{i,j}\right) .\left( \begin{array}{ccc}
k.\lambda _{n} &  & \\
 & \ddots & \\
 &  & k.\lambda _{n}
\end{array}\right) .\]
 So the result.\( \;  \)
\end{proof}
\begin{lem}
Up to a covering, and in a well-fitted coordinate chart, the brackets induced by \(
\Pi  \) can be written in a neighborhood of \( \Gamma  \) as:
\[
\begin{array}{cll}
\{x_{i},x_{j}\}&=&o_{2}(x)\\
\; \; \{\theta ,x_{i}\}&=&k(\theta )\lambda _{i}+o_{2}(x)\qquad \forall i=1\dots n,
\end{array}\]
where \( k\in C^{\infty }(S^{1}) \).
\end{lem}
\begin{proof}
For any \( \theta \in S^{1} \), \( H_{\theta } \) has distinct eigenvalues \(
\{k(\theta )\lambda _{i}\ ;\ i=1\dots n\} \) so \( k \) is smooth, and, if we denote \(
E_{\lambda _{i}} \) the associated characteristic spaces, they are supplementary
subvectorbundles of \( S^{1}\times \mathbb{R}^{n}\longmapsto S^{1} \), so that each
of them is either trivial, or diffeomorphic to the Mo\( \!  \)ebius band \( M \).

The covering
\[
\begin{array}{l}
C:\bf S^{1}\times \mathbb{R}^{n}\longmapsto S^{1}\times \mathbb{R}^{n}\\
\qquad \; (\widetilde{\theta },\widetilde{x})\; \; \longmapsto \; \; (2\widetilde{\theta },\widetilde{x}),
\end{array}\]
induces on the source space a Poisson structure \( \widetilde{\Pi } \) for which the
associated characteristic bundles \( E_{\widetilde{\lambda _{i}}} \) are all
trivial. As the topology of these subbundles is invariant by Poisson isomorphisms,
it is clear that two germs of Poisson structures \( \Pi  \) and \( \Pi ' \) on \(
S^{1}\times \mathbb{R}^{n} \) are isomorphic if and only each \( E_{\lambda _{i}} \)
is diffeomorphic to one of the \( E_{\lambda ^{'}_{i}} \)'s and \( \widetilde{\Pi }
\) is Poisson isomorphic to \( \widetilde{\Pi '} \).

Also, up to a covering, \( E_{\lambda _{i}} \) is trivial for any \( i \).
Then the choice of non vanishing smooth sections for each of them defines a
linear change of coordinates that diagonalizes the matrix \( H_{\theta } \).
\end{proof}
\begin{rem}
Of course, Poisson structures on \( S^{1}\times \mathbb{R}^{n} \) admitting non
trivial characteristic subbundles do exist, the basic exemple on \( S^{1}\times
\mathbb{R}^{2} \) is given by the following brackets:
\[
\begin{array}{l}
\{x_{1},x_{2}\}=\frac{x_{1}^{2}}{2}(\lambda cos^{2}\frac{\theta }{2}+\mu sin^{2}\frac{\theta }{2})+\frac{x^{2}_{2}}{2}(\lambda sin^{2}\frac{\theta }{2}+\mu cos^{2}\frac{\theta }{2})+x_{1}x_{2}(\mu -\lambda )cos\frac{\theta }{2}sin\frac{\theta }{2}\\
\; \{\theta ,x_{1}\}\; =x_{1}(\lambda cos^{2}\frac{\theta }{2}+\mu sin^{2}\frac{\theta }{2})+x_{2}(\mu -\lambda )cos\frac{\theta }{2}sin\frac{\theta }{2}\\
\; \{\theta ,x_{2}\}\; =x_{1}(\mu -\lambda )cos\frac{\theta }{2}sin\frac{\theta }{2}+x_{2}(\lambda sin^{2}\frac{\theta }{2}+\mu cos^{2}\frac{\theta }{2}).
\end{array}\]
Here, it is interesting to notice how these brackets simplify when \( \lambda =\mu =1 \).
\end{rem}
\begin{lem}
By re-parametrizing \( S^{1} \), one can put the brackets under the following form:
\[
\begin{array}{l}
\{x_{i},x_{j}\}=o_{2}(x)\\
\; \{\theta ,x_{i}\}\; =\mu _{i}x_{i}+o_{2}(x),\qquad  i=1\dots n,
\end{array}\ \]

\end{lem}
\begin{proof}
The formula $$ \chi (\theta )=\frac{2\pi }{\int _{0}^{2\pi }k(t)dt}\int _{0}^{\theta
}k^{-1}(t)dt $$  defines a diffeomorphism \( S^{1}\longmapsto S^{1} \) that gives
the announced re-parametrization.
\end{proof}
\begin{lem}
There exists a coordinate chart, defined on a neighborhood of \( \Gamma  \), in
which the brackets express as:
\[
\begin{array}{l}
\{x_{i},x_{j}\}=o_{2}(x)\\
\: \{\theta ,x_{i}\}\; =\mu _{i}x_{i},\qquad  i=1\dots n.
\end{array}\]
\end{lem}
\begin{proof}
Let \( X_\theta \) the hamiltonian vector field associated to the
(multi-)function \( \theta  \):
\[
X_\theta=\sum _{i=1\dots n}(\mu _{i}x_{i}+p_{i})\frac{\partial }{\partial x_{i}},\]
 where \( p_{i} \) are \( o_{2}(x) \). This vector field turns out to be linearisable, considered as a vector field in the variables $x$ with parameter $\theta$: there is a family of local diffeomorphisms $\phi_\theta ,$ depending smoothly on $\theta ,$ such that
\[
\phi _{\theta *}X_\theta=\sum _{i=1\dots n}\mu _{i}x_{i}\frac{\partial }{\partial x_{i}}\ \]
for every $\theta .$ Moreover we can choose $\phi_\theta $ 1-periodic in $\theta .$
So the change of coordinates $(x,\theta )\mapsto (\phi_\theta (x) , \theta)$ gives
the result.
\end{proof}
\begin{rem} Let us give an idea of how to prove that $X_{\theta}$ is linearisable:
first notice that a diffeomorphism $\phi_\theta=(\phi^1_\theta,\dots ,\phi^n_\theta)$ linearises $X_\theta$ if and only if 
 $$ {\mathcal L}_{X_\theta} \phi^i_\theta=\mu_i\phi^i_\theta \quad \forall i=1\dots  n, \forall \theta\in S^1.$$
 The first step consists in proving the result formally: write $X_\theta$ as 
 $$X_\theta=X^{(1)}_\theta+X^{(r)}_\theta+o_r({\bf x}),$$
 with $X^{(1)}_\theta, X^{(r)}_\theta$ and $o_r(x)$ denoting respectively the terms of order $1$, of order $r>1$ and of order higher than $r$ in the $x_i$'s variables, and show the existence of a formal diffeomorphism $\phi_\theta=Id+\sum_{i\in {\mathbb N}}\phi_\theta^{(r)}$ by induction on $r$, using non-resonnance conditions. Then, under Bruno's condition, one can show that the series $\phi_\theta$ converge, getting the result in the analytic context. In the smooth case, one can proceed as in $[RS]$ : what preceeds, leads us to the form $X_\theta=X_\theta^{(1)}+o_\infty(x)$ where $o_\infty(x)$ denotes a  vector field flat at $0\in{\mathbb R}^n$ for all $\theta\in S^1$, it is then possible to eliminate these flat terms by direct integration, flatness assuring the convergence of this integral.
\end{rem}
 
\begin{lem}
A last change of coordinates leads to the following brackets:
\[
\begin{array}{l}
\{x_{i},x_{j}\}=a_{i,j}x_{i}x_{j}+o_{\infty }(x),\qquad i,j=1\dots n\\
\; \{\theta ,x_{i}\}\; =\mu _{i}x_{i},\qquad \qquad \quad \, \; \qquad  i=1\dots n.
\end{array}\]
\end{lem}
\begin{proof}
Let \( u_{i,j}=\{x_{i},x_{j}\} \), Jacobi's identity \( \oint _{\theta
,x_{i,}x_{j}}\{\theta ,\{x_{i},x_{j}\}\}=0  \) is now just:
\[
(\mu _{i}+\mu _{j})u_{i,j}=X_{\theta }(u_{i,j}).\]
 First notice that \(
u_{i,j}^{0}=x_{i}x_{j} \) is a particular solution of this equation, and that \(
u_{i,j} \) being another solution, writing it as \( u_{i,j}=w_{i,j}x_{i}x_{j} \), \(
w_{i,j} \) has to be a first integral of \( X_{\theta } \).

Now, because of the dynamic of this vector field , the restriction \( (w_{i,j})_{\mid \bigcup E_{\mu _{i}}} \)only
depends on \( \theta  \) so \( w_{i,j}=k_{i,j}(\theta )+\tau _{i,j} \), where
\( \tau _{i,j} \) is flat on \( \Gamma  \) (notice that if all the \( \mu _{i} \)
's have the same sign, then \( \tau _{i,j} \) is necessarily zero because it
is a first integral of \( X_{\theta } \)).

Writing the terms of order one in Jacobi's identity \( \oint
_{x_{1},x_{i},x_{j}}\{x_{1},\{x_{i\, },x_{j}\}\}=0 \) with \( k_{1,i} \) and \(
k_{1,j} \) independent of \( \theta  \) implies that \( k_{i,j} \) is independent of
\( \theta  \) too, so we only have to find a change of coordinates making \( k_{1,i}
\) independent of \( \theta  \) for any \( i=1\dots n \). This diffeomorphism can be
defined by:
\[
(\theta ,x_{1},\dots ,x_{n})\longmapsto (\theta ,x_{1},\chi _{2}x_{2},\dots ,\chi _{n}x_{n}),\]
 where $$ \chi _{j}(\theta )=e^{\frac{1}{\lambda _{1}}\int ^{\theta }_{0}(k_{1,j}(t)-\overline{k}_{1,j})dt}
 $$
with \( \overline{k}_{1,j}=\frac{1}{2\pi }\int _{0}^{2\pi }k_{1,j}(t)dt \) for any
\( j=2\dots n \).
\end{proof}

\begin{rem} In the smooth context, when $n=2$, it is still possible to improve this result by eliminating flat terms, the proof is rather technical and specific to the case $n=2$, it can be found in $[Br]$.
\end{rem}
 
\begin{rem} In the case $n=1$ the normal form reduces to
$$\pi = cx\frac{\partial}{\partial \theta}\wedge \frac{\partial}{\partial x}$$
and we obtain, as a particular case, the form given in {[}R{]}.

In the case $n=2$ we can obtain a slightly better result: the (analytic and smooth) normal form is
$$\pi = \frac{\partial}{\partial \theta}\wedge (\mu_1x_1\frac{\partial}{\partial x_1}+
\mu_2x_2\frac{\partial}{\partial x_2})+ax_1x_2\frac{\partial}{\partial
x_1}\wedge\frac{\partial}{\partial x_2}$$ where $a$, and $\mu_1,\mu_2$ are real constants.
\end{rem}

\section{Geometric significance of the parameters.}

In that section we will see that the parameters $\mu_i$ and $a_{i,j},$ which appear
in the above normal form, are invariants of the Poisson structure (by this we mean
that two equivalent such Poisson structures must have normal forms with the same
$\mu_i$'s and $a_{i,j}$'s) and give them a geometrical interpretation.
The sum of the $\mu_i$'s gives the period $c$ of the modular vector field along the singular curve, the $n$-tuple $(\mu_1,\dots  ,\mu_n)$ caracterises, up to a multiplicative constant, the linear approximation of $\Pi$ along $\Gamma$, finally the $\mu_i$'s and $a_{i,j}$'s measure the (non Poisson) holonomy of the symplectic foliation when we move along this singular curve.

\noindent{$\bullet$ \textit{Period of the rotational and linear approximation along $\Gamma$}}

The modular vector field (see {[}W1{]}) \( D_{\omega }\Pi \) of \( \Pi \) with
respect to the volume form \( \omega =d\theta \wedge dx_{1}\wedge \dots \wedge dx_{n}
\) (coordinates of the normal form) has the following expression:
\[
D_{\omega }\Pi =(\sum _{i=1\dots n}\mu _{i}+o_{1}(x))\frac{\partial }{\partial \theta
}+o_{2}(x).\] As \( {\mathcal D}_{\omega }\Pi  \) is defined up to hamiltonian vector fields,
its restriction to \( \Gamma  \) depends only on \( \Pi  \), but \( {2\pi }/{\sum
_{i=1\dots n}\mu _{i}} \) appears to be the period of its flow on \( \Gamma  \). So, as
we have remarked in the introduction, it is an invariant of the Poisson structure.

We have seen that, at a purely linear level, the set \( \{\mu _{1},\dots \, ,\mu _{n}\}
\) was defined up to a multiplicative constant. Because the linear part of the
Poisson structure along the singular curve is (up isomorphism) an invariant this
set, up to multiplication by a constant, is also an invariant. If we use the
preceding result we obtain that the $\mu_i$ themselves are invariants. We recall
that, in the local study ({[}DZ{]}) the $\mu_i$ are only invariant up to
multiplication by a constant. Also this local study shows that  the $\mu_i$ (up to
multiplication by a constant) measure the way the symplectic leaves behave on a
transversal to the singular curve: the traces of these leaves on such a transversal
give a figure diffeomorphic to the phase portrait of the linear vector field $\sum
_{i=1\dots n}\mu_{i}x_{i}{\partial }/{\partial x_{i}}.$

\noindent {$\bullet$ \textit{The symplectic foliation}}

Before showing how the $a_{i,j}$'s are involved in the symplectic foliation of $\Pi$, let us show that they are Poisson invariants : consider two Poisson
structures \( \Pi  \) and \( \Pi ' \) on a neighborhood of \( \Gamma  \), and \(
\Phi  \) a Poisson isomorphism from \( \Pi  \) to \( \Pi ' \)~. Put \( \Pi  \) and
\( \Pi ' \) under normal form, as the \( \mu _{i} \) 's are Poisson invariants, we
can suppose that \( \mu _{i}=\mu '_{i} \) for any \( i=1\dots n \). Then, writing down
the equations for \( \Phi  \) to be Poisson, we get, with evident notations, \(
\forall i,j=1...n \):
\[
\begin{array}{c}
\sum ^{n}_{k=1}\mu _{k}x_{k}(\frac{\partial \Phi _{0}}{\partial \theta }\frac{\partial \Phi _{i}}{\partial x_{k}}-\frac{\partial \Phi _{0}}{\partial x_{k}}\frac{\partial \Phi _{i}}{\partial \theta })+\sum ^{n}_{k,l=1}a_{k,l}x_{k}x_{l}(\frac{\partial \Phi _{0}}{\partial x_{k}}\frac{\partial \Phi _{i}}{\partial x_{l}})=\mu _{i}\Phi _{i}\qquad (I)_{i}\\
\quad \sum _{k=1}^{n}\mu _{k}x_{k}(\frac{\partial \Phi _{i}}{\partial \theta
}\frac{\partial \Phi _{j}}{\partial x_{k}}-\frac{\partial \Phi _{i}}{\partial
x_{k}}\frac{\partial \Phi _{j}}{\partial \theta })+\sum
^{n}_{k,l=1}a_{k,l}x_{k}x_{l}\frac{\partial \Phi }{\partial x_{k}}\frac{\partial
\Phi }{\partial x_{l}}=a'_{i,j}\Phi _{i}\Phi _{j}\quad (II)_{i,j}.
\end{array}\]
Differentiating equation \( (I)_{i} \) with respect to \( \theta  \), and evaluating
it onto \( \Gamma  \), we obtain:
\[
\frac{\partial \Phi _{i}}{\partial \theta }(\theta ,0,0)=0\qquad \forall \theta \in
S^{1}.\]
 Differentiating it with respect to \( x_{j} \), and evaluating it onto \(
\Gamma  \), we get:
\[
\mu _{j}(\frac{\partial \Phi _{0}}{\partial \theta }\frac{\partial \Phi
_{i}}{\partial x_{j}}-\frac{\partial \Phi _{0}}{\partial x_{j}}\frac{\partial \Phi
_{j}}{\partial \theta })(\theta ,0,0)=\mu _{i}\frac{\partial \Phi _{i}}{\partial
x_{j}}(\theta ,0,0)\qquad \forall \theta \in S^{1}.\]

But \( \Phi  \)\( \mid _{\Gamma } \) must exchange modular vector fields, so \( \Phi
_{1}(\theta ,0,0)=\theta +constant \) and this equation becomes:
\[
\frac{\partial \Phi _{i}}{\partial x_{j}}(\theta ,0,0)=0\qquad \forall \theta \in
S^{1}.\]

In conclusion, we can write \( \Phi  \) as:

\[
\begin{array}{c}
\Phi _{0}=\theta +o_{2}(x)\qquad \qquad \qquad \qquad \quad \\
\Phi _{i}=u_{i}(\theta )x_{i}+o_{2}(x)\qquad \forall i=1...n
\end{array}\]

And in \( (II)_{i,j} \), the terms of order two in the $x_i$'s variables give the relation:
\[
\frac{\partial }{\partial \theta }(\mu _{j}\, ln\mid \! u_{i}\! \mid +\mu _{i}\,
ln\mid \! u_{j}\! \mid )=a_{i,j}-a'_{i,j}\]
 Here, the left side is the derivative of a function defined on \( S^{1} \),
so its integral over \( S^{1} \) is zero, and \( a_{i,j}=a'_{i,j} \).

Let us now describe the symplectic foliation of $\Pi$, from on (in the smooth case), it will be supposed that {\it no flat terms appear in the normal form}, in other words, we suppose that
$$ \Pi={\partial \theta}  \wedge \Bigl( \sum_{i=1\dots  n} \mu_i x_i {\partial x_i} \Bigr) + \sum_{1\leq i < j\leq n}a_{i,j}x_i x_j {\partial x_i} \wedge {\partial x_j}.$$

We will only describe the symplectic foliation on the open set $P^+:=\{(\theta,x)\in S^1\times {\mathbb R}^n/x_i > 0, \forall i=1 \dots  n\}$, however, the following reasonning allows to describe it entirely: indeed if, for any subset $I\subset \{1,\dots  ,n\}$ we denote $P_I=\{(\theta,x)\in S^1\times {\mathbb R}^n/x_i=0, \forall i\notin I\}$, one easily sees that $P_I$ is a Poisson submanifold of the presently studied type (that contains itself Poisson sumbmanifolds $P_J$ with $J \subset I$), so we can apply what follows to $P_I^+:=\{(\theta,x)\in S^1\times {\mathbb R}^n/x_i=0, \forall i\notin I, x_i>0\quad else\}$. Moreover, it will appear clearly that $P^+$ is one of the connected componant of the regular open set of $\Pi$ $-$all these component being isomorphic as the reflexions $(\theta,x_1,\dots ,x_i,\dots  ,x_n) \longmapsto (\theta,x_1,\dots ,-x_i,\dots ,x_n)$ are Poisson$-$ and that the complement of this regular set is precisely composed by the $P_I$ for $I \subset \{1 \dots  n\}, I \neq \{1 \dots  n\}$.

To put the foliation induced by $\Pi$ on $P^+$ in evidence, one first apply the diffeomorphism 
$$ \begin{array}{cccc}
     L: & P^+ & \longmapsto & S^1 \times {\mathbb R} ^n \\
          &(\theta,x_1,\dots ,x_n) & \longmapsto & (\theta,\ln x_1,\dots ,\ln x_n), \\
    \end{array} $$
as 
 $$ \{\theta,\ln x_i\}=\mu_i
    \{\ln x_i,\ln x_j\}=a_{i,j}, $$
in the induced coordinates (that we shall denote by $(\theta,{\bf \overline{x}})=(\theta,\overline{x}_1,\dots ,\overline{x}_n)=(\theta, \ln x_1,\dots ,\ln x_n)$), the matrix $( \Pi_{i,j} )$ of the brackets has constant coefficents:
$$ ( \Pi_{i,j} )=
 \left(
         \begin{array}{cccc}
                                0      &   -\mu_1 & \dots  & -\mu_n     \\
                                \mu_1                              \\
                                \vdots &       &  (a_{i,j})     \\
                                \mu_n                              \\
          \end{array}
       \right). $$
Denote $\mu$ the vector $(\mu_1,\dots  ,\mu_n)$ and $2s$ the rank of the matrix $(a_{i,j})$, the theorem of spectral decomposition for anti-autoadjoint operators assures the existence of a linear isomorphism $\phi:{\mathbb R}^n \longmapsto {\mathbb R}^n$, whose matrix will be denoted $(\phi_{i,j} )$ that conjugates $(a_{i,j})$ to a matrix of the following type 
$$\left(
\begin{array}{cccc}
 \left( \begin{array}{cc}
          0 & -1 \\
          1 & 0 \end {array}
          \right)           &      &    O                       &    O    \\
                            &\ddots&                            &  \vdots \\
              O             &      &  \left( \begin{array}{cc}   
                                              0 & -1 \\
                                              1 & 0 \end {array}
                                                  \right)       & O          \\
              O             &\dots &              O               & \left( \begin{array}{ccc}   
                                                                                                                                0     & \dots  &0     \\
                                                                                                                                \vdots&       &\vdots \\
                                                                                                                                0     & \dots  &0 \end {array}
                                                                                                                               \right)                            \\
                                       
  \end{array}                                     
  \right) . $$
Moreover, if $\mu$ belongs to \`a $Im(a_{i,j})$, one can choose $\mu$ as an element of the basis for which $(a_{i,j})$ takes that form, this incites us to distinguish two cases :

\underbar{$1^{st}$ case}: if $\mu \in Im(a_{i,j})$

Chosing $\phi$ as described below, the diffeomorphism $(\theta,{\bf x})\longmapsto(\theta,\phi({\bf x}))$ conjugates $(\Pi_{i,j})$ to the following matrix:
$$\left(
         \begin{array}{ccccr}
                     0               & \begin{array}{cc} 
                                                0 & -1 \end{array}        &         &         \dots                                    &     0 \mbox{ }            \\
                    \begin{array}{c}       
                     0   \\   
                     1   \\    
                    \end{array}      & \left( \begin{array}{cc}   
                                                 0 & -1 \\
                                                 1 & 0 \end {array}
                                       \right)                            &         &    O                                                                        \\
                                     &                                    & \ddots  &                                                &     \vdots \mbox{ }        \\
                     \vdots          &                  O                 &         &  \left( \begin{array}{cc}   
                                                                                                0 & -1 \\
                                                                                                1 & 0 \end {array}
                                                                                       \right)                                                                    \\           
                    \begin{array}{c} \\ \\
                     0 \end{array}    &                                   &\begin{array}{c} \\ \\
                                                                              \dots  \end{array} &  \begin{array}{cc} \\ \\
                                                                                                & 0 \end{array}                       & \left( \begin{array}{ccc}   
                                                                                                                                0     & \dots  &0     \\
                                                                                                                                \vdots&       &\vdots \\
                                                                                                                                0     & \dots  &0 \end {array}
                                                                                                                               \right)                            \\
                                       
  \end{array}                                     
  \right) . $$
Here, the null matrix, in the lower right corner is of order $n-2s$.

Denote $(\theta,q_1,p_1,\dots ,q_s,p_s,\tilde{x}_{2s+1},\dots ,\tilde{x}_{n})$ the coordinates induced by $\phi \circ L$, the foliation is then spanned by $ {\partial \theta}+{\partial q_1},
{\partial p_1}, {\partial q_2}, {\partial p_2}, \dots 
  {\partial q_s},  {\partial p_s}$, so that this foliation is not parallel to $S^1\times\{0\}$, letting a holonomy phenomenon appear along $\Gamma$. The set of leaves can be parametrised by $(q_1,\tilde{x}_{2s+1},\dots ,\tilde{x}_n)\in[0,2\pi[\times {\mathbb R}^{n-2s}$.
 
Moreover, if $\mathcal{F}_{(q_1,\tilde{x}_{2s+1},\dots ,\tilde{x}_n)}$ denotes the leaf through the point  $(0,q_1,0,\dots ,0,\tilde{x}_{2s+1},\dots ,\tilde{x}_n)\in S^1\times {\mathbb R} ^n$,
 then $\mathcal{F}_{(q_1,\hat{x}_{2s+1},\dots ,\hat{x}_n)}$ is parametrised by $(t_1,\dots ,t_{2s})$: 
 $$\mathcal{F}_{(q_1,\tilde{x}_{2s+1},\dots ,\tilde{x}_n)}=\{(t_1,q_1+t_1,t_2,t_3,\dots ,t_{2s},\tilde{x}_{2s+1},\dots ,\tilde{x}_n)\in S^1\times {\mathbb R} ^n /t_1\dots  t_{2s}\in{\mathbb R}\},$$

 In the coordinates of the normal form, it appears that the leaf $F_{(x_1,\dots ,x_n)}$ through the point $(0,x_1,\dots ,x_n)$ admits the following  parametrisation:
 $$F_{(x_1,\dots ,x_n)}=\{(t_1,(x_ie^{\sum_{j=1\dots  2s} \psi_{i,j}t_j})_{i=1\dots  n}) \in P^+ /t_1\dots  t_{2s} \in{\mathbb R}\},$$
 where $(\psi_{i,j})$ denotes the inverse matrix of $(\phi_{i,j})$.
  
 \underbar{$2^{nd}$ case}: if $\mu \notin Im(a_{i,j})$
 
 We follow the same reasoning: $L$ d\'efinit a chart on $P^+$ in which the matrix $(\Pi_{i,j})$ 
 can be conjugated by a linear isomorphism $\phi$ to the following
 $$\left(
         \begin{array}{ccccr}
                     0               & \begin{array}{cc} 
                                                 0 & 0 \end{array}        &         &         \dots                                    & \begin{array}{ccc} -1 & \dots  &  \end{array}\\
                    \begin{array}{c}       
                     \vdots \\   
                            \\    
                    \end{array}      & \left( \begin{array}{cc}   
                                                 0 & -1 \\
                                                 1 & 0 \end {array}
                                       \right)                            &         &    O                                                                        \\
                                     &                                    & \ddots  &                                                &     \vdots \mbox{ }        \\
                     \vdots          &                  O                 &         &  \left( \begin{array}{cc}   
                                                                                                0 & -1 \\
                                                                                                1 & 0 \end {array}
                                                                                       \right)                                       &                            \\           
                    \begin{array}{c} 1\\ \vdots \\
                     0 \end{array}    &                                   &\begin{array}{c} \\ \\
                                                                              \dots  \end{array} &  \begin{array}{cc} \\ \\
                                                                                                & 0 \end{array}                       & \left( \begin{array}{ccc}   
                                                                                                                                0     & \dots  &0     \\
                                                                                                                                \vdots&       &\vdots \\
                                                                                                                                0     & \dots  &0 \end {array}
                                                                                                                               \right)                            \\
                                       
  \end{array}                                     
  \right) . $$
  In the coordinates induced by $\phi \circ L$ , that we shall denote by $(\theta,q_1,p_1,\dots ,q_s,p_s,\tilde{x}_{2s+1},\dots ,\tilde{x}_{n})$, the foliation is spanned by  $\partial \theta, {\partial q_1},
{\partial p_1}, {\partial q_2}, {\partial p_2}, \dots ,
  {\partial q_s},  {\partial p_s}, \partial \tilde{x}_{2s+1}$, so it is parallel to $S^1\times\{0\}$ and there is no holonomy in that case.
  
  Let us denote $\mathcal{F}_{(\tilde{x}_{2s+2},\dots ,\tilde{x}_n)}$ the leaf through the point  $(0,0,\dots ,0,\tilde{x}_{2s+2},\dots ,\tilde{x}_n)\in S^1\times {\mathbb R} ^n$, we get:
$$\mathcal{F}_{(\tilde{x}_{2s+2},\dots ,\tilde{x}_n)}=\{(e^{it_{0}},t_1,\dots ,t_{2s},\tilde{x}_{2s+1}=t_{2s+1},\tilde{x}_{2s+2},\dots ,\tilde{x}_n) \in P^+/t_0,\dots ,t_{2s+1} \in {\mathbb R}\}$$
   In the coordinates of the normal form, denote $F_{(x_1,\dots ,x_n)}$ the leaf through the point $(0,x_1,\dots ,x_n)$, it admits a parametrisation:
 $$F_{(x_1,\dots ,x_n)}=\{(t_0,(x_ie^{\sum_{j=1\dots  2s+1} \psi_{i,j}t_j})_{i=1\dots  n}) \in P^+ /t_0\dots  t_{2s+1} \in{\mathbb R}\},$$
 where $(\psi_{i,j})$ denotes the inverse matrix of $(\phi_{i,j})$.

\end{document}